\documentclass[10pt,sumlimits,namelimits]{article}
\newcommand{\nc}{\newcommand}
\newcommand{\rnc}{\renewcommand}
\usepackage{amsthm}
\usepackage{amsmath}
\usepackage{amsfonts}
\usepackage{amssymb}
\usepackage{color}
\usepackage[colorlinks=true,urlcolor=blue]{hyperref}
\usepackage[applemac]{inputenc}

\nc{\N}{\mathbb{N}}
\nc{\Z}{\mathbb{Z}}
\nc{\D}{\mathbb{D}}
\nc{\Q}{\mathbb{Q}}
\nc{\R}{\mathbb{R}}
\nc{\C}{\mathbb{C}}
\nc{\vS}{{\cal S}}
\nc{\vphi}{\varphi}
\nc{\eps}{\varepsilon}
\nc{\dsp}{\displaystyle}
\nc{\ovl}{\overline}
\nc{\udl}{\underline}
\nc{\vlim}{\lim\limits}
\nc{\vlimsup}{\limsup\limits}
\nc{\vliminf}{\liminf\limits}
\nc{\vsup}{\sup\limits}
\nc{\vinf}{\inf\limits}
\nc{\vint}{\int\limits}
\nc{\inj}{\hookrightarrow}
\nc{\tends}{\longrightarrow}
\nc{\weak}{\rightharpoonup}
\nc{\w}{{\textsl w}}
\nc{\loc}{{\rm loc}}
\nc{\rad}{{\rm rad}}
\rnc{\le}{\leqslant}
\rnc{\ge}{\geqslant}
\rnc{\Re}{{\rm Re}}
\rnc{\Im}{{\rm Im}}

\numberwithin{equation}{section}
\rnc{\theequation}{\thesection.\arabic{equation}}

\newtheorem{thm}{Theorem}[section]
\newtheorem{prop}[thm]{Proposition}

\newtheorem{lem}[thm]{Lemma}

\theoremstyle{definition}
\newtheorem{rmk}[thm]{Remark}

\newenvironment{proof*}{\noindent{\bf Proof.}}{\qed}
\newenvironment{vproof}[1]{\noindent{\bf Proof #1}}{\qed}

\title{\Huge \sc Necessary Conditions and Sufficient Conditions for Global Existence in the Nonlinear Schrödinger Equation}
\author{\sc Pascal Bégout}
\date{}

\begin{document}

\maketitle

\begin{center}
Laboratoire Jacques-Louis Lions \\
Université Pierre et Marie Curie \\
Boîte Courrier 187 \\
4, place Jussieu 75252 Paris Cedex 05, FRANCE \bigskip \\
{\footnotesize e-mail\:: }\htmladdnormallink{{\footnotesize\udl{\tt{begout@ann.jussieu.fr}}}}
{mailto:begout@ann.jussieu.fr}
\end{center}

\begin{abstract}
In this paper, we consider the nonlinear Schrödinger equation with the super critical power of nonlinearity in the attractive case. We give a sufficient condition and a necessary condition to obtain global or blowing up solutions. These conditions coincide in the critical case, thereby extending the results of Weinstein \cite{MR84d:35140,MR87i:35026}. Furthermore, we improve a blow-up condition.
\end{abstract}

\baselineskip .7cm

\section{Introduction and notations}
\label{introduction}

{\let\thefootnote\relax\footnotetext{2000 Mathematics Subject Classification: 35Q55}}

We consider the following nonlinear Schrödinger equation,
\begin{gather}
 \left\{
  \begin{split}
   \label{nls}
    i\frac{\partial u}{\partial t}+\Delta u+\lambda|u|^\alpha u & = 0,\; (t,x)\in(-T_*,T^*)\times \R^N, \\
                                                           u(0) & = \vphi, \mbox{ in }\R^N,
  \end{split}
 \right.
\end{gather}
where $\lambda\in \R,$ $0\le\alpha<\dfrac{4}{N-2}$ $(0\le\alpha<\infty$ if $N=1)$ and $\vphi$ a given initial data.

It is well-known that for every $\vphi\in H^1(\R^N),$ (\ref{nls}) has a unique solution $u\in C((-T_*,T^*);H^1(\R^N))$ which 
satisfies the blow-up alternative and the conservation of charge and energy. In other words, if $T^*<\infty$ then
$\vlim_{t\nearrow T^*}\|u(t)\|_{H^1}=\infty.$ In the same way, if $T_*<\infty$ then $\vlim_{t\searrow-T_*}\|u(t)\|_{H^1}
=\infty.$ And for all $t\in(-T_*,T^*),$ $\|u(t)\|_{L^2}=\|\vphi\|_{L^2}$ and $E(u(t))=E(\vphi),$ where $E(\vphi)\stackrel{\rm
def}{=}\frac{1}{2}\|\nabla\vphi\|_{L^2}^2-\frac{\lambda}{\alpha+2}\|\vphi\|_{L^{\alpha+2}} ^{\alpha+2}.$ If $\vphi\in
X\stackrel{\rm def}{=}H^1(\R^N)\cap L^2(|x|^2;dx)$ then $u\in C((-T_*,T^*);X).$ Moreover, if $\lambda\le0,$ if
$\alpha<\frac{4}{N}$ or if $\|\vphi\|_{H^1}$ is small enough then $T^*=T_*=\infty$ and
$\|u\|_{L^\infty(\R;H^1)}<\infty.$ Finally, is $\alpha\ge\frac{4}{N}$ then there exist initial values $\vphi\in H^1(\R^N)$ such
that the corresponding solution of (\ref{nls}) blows up in finite time. See Cazenave \cite{caz1}, Ginibre and Velo
\cite{MR82c:35059,MR82c:35057,MR82c:35058,MR87b:35150}, Glassey \cite{MR57:842}, Kato \cite{MR88f:35133}.

In the attractive and critical case $(\lambda>0$ and $\alpha=\frac{4}{N}),$ there is a sharp condition to obtain global solutions
(see Weinstein \cite{MR84d:35140,MR87i:35026}). It is given in terms of the solution of a related elliptic problem. But in the
super critical case $(\alpha>\frac{4}{N}),$ we only know that there exists $\eps>0$ sufficiently small such that if
$\|\vphi\|_{H^1}\le\eps,$ then the corresponding solution is global in time.

In this paper, we try to extend the results of Weinstein \cite{MR84d:35140,MR87i:35026} to the super critical case
$\alpha>\frac{4}{N}.$ As we will see, we are not able to establish such a result, but we can give two explicit real values
functions $\gamma_*$ and $r_*$ with $0<\gamma_*<r_*$ such that if $\|\vphi\|_{L^2}\le\gamma_*(\|\nabla\vphi\|_{L^2}),$ then the
corresponding solution is global in time. Furthermore, for every $(a,b)\in(0,\infty)\times(0,\infty)$ such that $a>r_*(b),$ there
exists $\vphi_{a,b}\in H^1(\R^N)$ with $\|\vphi_{a,b}\|_{L^2}=a$ and $\|\nabla \vphi_{a,b}\|_{L^2}=b$ such that the associated
solution blows up in finite time for both $t<0$ et $t>0$ (see Theorem \ref{thmexplosion} below). Despite of the fact we do not
obtain a sharp condition (since $\gamma_*<r_*),$ we recover the results of Weinstein \cite{MR84d:35140,MR87i:35026} as
$\alpha\searrow\frac{4}{N}.$ Setting ${\cal A}=\{\vphi\in H^1(\R^N);$ $\|\vphi\|_{L^2}\le\gamma_*(\|\nabla\vphi\|_{L^2})\},$ it
follows that for every $\vphi\in{\cal A},$ the corresponding solution of (\ref{nls}) is global in time and uniformly bounded in
$H^1(\R^N).$ It is interesting to note that ${\cal A}$ is an unbounded subset of $H^1(\R^N)$ as for the case
$\alpha=\frac{4}{N}.$ We also improve some results about blow-up (Theorems \ref{explosionX} and \ref{explosionXrad}).

This paper is organized as follows. In Section \ref{blowup}, we give a sufficient blow-up condition. In Section
\ref{sharpestimate}, we recall the best constant in a Gagliardo-Nirenberg's inequality. In Section \ref{cns}, we give the main
result of this paper, that is necessary conditions and sufficient conditions to obtain global solutions. In Section
\ref{proofthmexplosion}, we prove the result given in Section \ref{cns}.

The following notations will be used throughout this paper.
$\Delta=\sum\limits_{j=1}^N\frac{\partial^2}{\partial x_j^2}$ and we denote  by $B(0,R),$ for $R>0,$ the ball of $\R^N$ of
center $0$ with radius $R.$
For $1\le p\le\infty,$ we design by $L^p(\R^N) = L^p(\R^N;\C),$ with norm $\|\: .\:\|_{L^p},$ the usual Lebesgue spaces and
by $H^1(\R^N) = H^1(\R^N;\C),$ with norm $\|\: .\:\|_{H^1},$ the Sobolev space.
 For $k\in\N\cup\{0\}$ and $0<\gamma<1,$ we
denote by
$C^{k,\gamma}(\R^N)=C^{k,\gamma}(\R^N;\C)$ the Hölder spaces and we introduce the Hilbert space $X=\left\{\psi\in H^1(\R^N;\C);\;
\|\psi\|_X<\infty\right\}$ with norm
$\|\psi\|_X^2=\|\psi\|_{H^1(\R^N)}^2+\vint_{\R^N}|x|^2|\psi(x)|^2dx.$
For a normed functional space $E\subset L^1_\loc(\R^N;\C),$ we denote by $E_\rad$ the space of functions $f\in E$ such that
$f$ is spherically symmetric. $E_\rad$ is endowed with the norm of $E.$
Finally, $C$ are auxiliary positive constants.

\section{Blow-up}
\label{blowup}

The first two results are an improvement of a blow-up condition (see Glassey \cite{MR57:842}, Ogawa and Tsutsumi
\cite{MR92k:35262}). We know that if a solution has a negative energy, then it blows up in finite time. We extend this result for
any nontrivial solution with nonpositive energy.

\begin{thm}
\label{explosionX}
Let $\lambda >0,$ $\dfrac{4}{N}<\alpha<\dfrac{4}{N-2}$ $(4<\alpha<\infty$ if $N=1)$ and $\vphi\in X,$ $\vphi\not\equiv0.$ If
$E(\vphi)\le0$ then the corresponding solution $u\in C((-T_*,T^*);X)$ of $(\ref{nls})$ blows up in finite time for both $t>0$ and
$t<0.$ In other words, $T^*<\infty$ and $T_*<\infty.$
\end{thm}

\begin{thm}
\label{explosionXrad}
Let $\lambda >0,$ $N\ge2,$ $\dfrac{4}{N}<\alpha<\dfrac{4}{N-2}$ $(2<\alpha\le4$ if $N=2)$ and $\vphi\in H^1_\rad(\R^N),$
$\vphi\not\equiv0.$ If $E(\vphi)\le0$ then the corresponding solution $u\in C((-T_*,T^*);H^1(\R^N))$ of $(\ref{nls})$ blows up in finite time for both $t>0$ and $t<0.$ In other words, $T^*<\infty$ and $T_*<\infty.$
\end{thm}

\begin{rmk}
When $E(\vphi)=0,$ the conclusion of Theorems \ref{explosionX} and \ref{explosionXrad} is false for $\alpha=\dfrac{4}{N}.$
Indeed, let $\vphi\in X_\rad,$ $\vphi\not\equiv0,$ be a solution of $-\Delta \vphi+\vphi=\lambda|\vphi|^{\frac{4}{N}}\vphi,$ in
$\R^N.$ Then $E(\vphi)=0$  from (\ref{groundstate3}) but $u(t,x)=\vphi(x)e^{it}$ is the solution of (\ref{nls}) and so
$T^*=T_*=\infty.$
\end{rmk}

Similar results exist for the critical case. See Nawa \cite{nawa,MR98m:35192}. It is shown that if $\vphi\in H^1(\R^N)$
satisfies $E(\vphi)<\dfrac{(\vphi',i\vphi)^2}{\|\vphi\|^2_{L^2}},$ when $N=1,$ or if $E(\vphi)<0,$ when $N\ge2,$ then the
corresponding solution of (\ref{nls}) blows up in finite time or grows up at infinity, the first case always occurring when $N=1.$
Here, $( \: , \: )$ denotes the scalar product in $L^2(\R^N).$ See also Nawa \cite{MR95g:35195,MR99m:35235}. Note that in the
case $N=1,$ the result of Nawa \cite{MR98m:35192} slightly improves that of Ogawa and Tsutsumi \cite{MR91f:35026}, since it
allows to make blow-up some solution with nonnegative energy.
\medskip

\begin{vproof}{of Theorem \ref{explosionX}.}
We argue by contradiction. Set for every $t\in(-T_*,T^*),$ $h(t)=\|xu(t)\|_{L^2}^2.$ Then $h\in C^2((-T_*,T^*);\R)$ and
\begin{gather}
\label{explosionX1}
\forall t\in(-T_*,T^*),\; h''(t)=4N\alpha E(\vphi)-2(N\alpha-4)\|\nabla u(t)\|^2_{L^2}
\end{gather}
(Glassey \cite{MR57:842}). Since $E(\vphi)\le0,$ we have by Gagliardo-Nirenberg's inequality (Proposition \ref{constGN}) and
conservation of energy and charge,
$
\|\nabla u(t)\|^2_{L^2}\le\frac{2\lambda}{\alpha+2}\|u(t)\|_{L^{\alpha+2}}^{\alpha+2}\le C\|\nabla
u(t)\|^{\frac{N\alpha}{2}}_{L^2},
$
for every $t\in(-T_*,T^*).$ Since $\alpha>\dfrac{4}{N}$ and $\vphi\not\equiv0,$ we deduce that
$\vinf_{t\in(-T_*,T^*)}\|\nabla u(t)\|_{L^2}>0$ and with (\ref{explosionX1}), we obtain
$$
\forall t\in(-T_*,T^*),\; h''(t)\le-C.
$$
So, if $T^*=\infty$ or if $T_*=\infty$ then there exists $S\in(-T_*,T^*)$ with $|S|$ large enough such that $h(S)<0$ which is absurd since $h>0.$ Hence the result.
\medskip
\end{vproof}
\\
\begin{vproof}{of Theorem \ref{explosionXrad}.}
For $\Psi\in W^{4,\infty}(\R^N;\R),$ $\Psi\ge0,$ we set
$$
\forall t\in(-T_*,T^*),\; V(t)=\vint_{\R^N}\Psi(x)|u(t,x)|^2dx.
$$
We know that there exists $\Psi\in W^{4,\infty}(\R^N;\R),$ $\Psi\ge0,$ such that $V\in C^2((-T_*,T^*);\R)$ and
$$
\forall t\in(-T_*,T^*),\; V''(t)\le 2N\alpha E(\vphi)-2(N\alpha-4)\|\nabla u(t)\|_{L^2}^2,
$$
(see the proof of Theorem 2.7 of Cazenave \cite{caz2} and Remark 2.13 of this reference). We conclude in the same way that for
Theorem \ref{explosionX}.
\end{vproof}

\section{Sharp estimate}
\label{sharpestimate}

In this section, we recall the sharp estimate in a Gagliardo-Nirenberg's inequality (Proposition \ref{constGN}) and a
result concerning the {\it ground states}.

Let $\lambda>0,$ $\omega>0$ and $0<\alpha<\dfrac{4}{N-2}$ $(0<\alpha<\infty$ if $N=1).$ We consider the following elliptic
equations.
\begin{gather}
\label{groundstate}
 \begin{cases}
  -\Delta R+R=|R|^\alpha R, \mbox{ in } \R^N, \\
  R\in H^1(\R^N;\R),\; R\not\equiv0,
 \end{cases}
\end{gather}
\begin{gather}
\label{groundstate'}
 \begin{cases}
  -\Delta \Phi+\omega\Phi=\lambda|\Phi|^\alpha\Phi, \mbox{ in } \R^N, \\
  \Phi\in H^1(\R^N;\R),\; \Phi\not\equiv0.
 \end{cases}
\end{gather}

It is well-known that the equation (\ref{groundstate'}) possesses at less one solution $\psi.$ Furthermore, each solution $\psi$
of (\ref{groundstate'}) satisfies $\psi\in C^{2,\gamma}(\R^N)\cap W^{3,p}(\R^N),$ $\forall\gamma\in(0,1),$ $\forall
p\in[2,\infty),$ $|\psi(x)|\le Ce^{-\delta|x|},$ for all $x\in\R^N,$ where $C$ and $\delta$ are two positive constants which do
not depend on $x,$ $\vlim_{|x|\to\infty}|D^\beta\psi(x)|=0,$ $\forall|\beta|\le2$ multi-index. Finally, $\psi$ satisfies the
following identities.
\begin{align}
\label{groundstate1}
            \|\nabla\psi\|_{L^2}^2 & =\dfrac{\omega N\alpha}{4-\alpha(N-2)}\|\psi\|_{L^2}^2, \medskip \\
\label{groundstate2}
\|\psi\|_{L^{\alpha+2}}^{\alpha+2} & =\dfrac{2\omega(\alpha+2)}{\lambda(4-\alpha(N-2))}\|\psi\|_{L^2}^2, \medskip \\
\label{groundstate3}
\|\psi\|_{L^{\alpha+2}}^{\alpha+2} & =\dfrac{2(\alpha+2)}{\lambda N\alpha}\|\nabla\psi\|_{L^2}^2. 
\end{align}
Such solutions are called {\it bound states} solutions. Furthermore, (\ref{groundstate'}) has a unique solution $\Phi$ satisfying
the following additional properties. $\Phi\in\vS_\rad(\R^N;\R);$ $\Phi>0$ over $\R^N;$ $\Phi$ is decreasing with respect to
$r=|x|;$ for every multi-index $\beta\in\N^N,$ there exist two constants $C>0$ and $\delta>0$ such that for every $x\in\R^N,$
$|\Phi(x)|+|D^\beta\Phi(x)|\le Ce^{-\delta |x|}.$ Finally, for every solution $\psi$ of (\ref{groundstate'}), we have
\begin{gather}
\label{groundstate4}
\|\Phi\|_{L^2}\le\|\psi\|_{L^2}.
\end{gather}
Such a solution is called a {\it ground state} of the equation (\ref{groundstate'}).

Equation (\ref{groundstate'}) is studied in the following references. Berestycki, Gallouët and Kavian~\cite{MR85e:35041};\linebreak Berestycki and Lions \cite{MR84h:35054a,MR84h:35054b} ; Berestycki, Lions and Peletier \cite{MR83e:35009}; Gidas, Ni and Nirenberg \cite{MR84a:35083}; Jones and Küpper \cite{MR87h:35261}; Kwong \cite{MR90d:35015}; Strauss \cite{MR56:12616}. See also Cazenave \cite{caz1}, Section 8.

\begin{prop}
\label{constGN}
Let $0<\alpha<\dfrac{4}{N-2}$ $(0<\alpha<\infty$ if $N=1)$ and $R$ be the ground state solution of $(\ref{groundstate}).$
Then the best constant $C_*>0$ in the Gagliardo-Nirenberg's inequality,
\begin{gather}
\label{constGN1}
\forall f\in H^1(\R^N),\; \|f\|_{L^{\alpha+2}}^{\alpha+2}\le C_*\|f\|_{L^2}^\frac{4-\alpha(N-2)}{2}\|\nabla
f\|_{L^2}^{\frac{N\alpha}{2}},
\end{gather}
is given by
\begin{gather}
\label{constGN2}
C_*=\frac{2(\alpha+2)}{N\alpha}\left(\frac{4-\alpha(N-2)}{N\alpha}\right)^{\frac{N\alpha-4}{4}}\|R\|_{L^2}^{-\alpha}.
\end{gather}
\end{prop}

See Weinstein \cite{MR84d:35140} for the proof in the case $N\ge2.$ See also Lemma 3.4 of Cazenave \cite{caz2} in the case $\alpha=\dfrac{4}{N}.$ But for convenience, we give the proof. It makes use of a compactness result which is an adaptation of the compactness lemma due to Strauss (Strauss \cite{MR56:12616}).

\begin{vproof}{of Proposition \ref{constGN}.}
We define for every $f\in H^1(\R^N),$ $f\not\equiv0,$ the functional
$$
J(f)=\dfrac{\|f\|_{L^2}^{\frac{4-\alpha(N-2)}{2}}\|\nabla f\|_{L^2}^{\frac{N\alpha}{2}}}{\|f\|_{L^{\alpha+2}}^{\alpha+2}},
$$
and we set $\sigma=\vinf_{f\in H^1\setminus\{0\}}J(f)$. Then $\sigma\in(0,\infty)$ by (\ref{constGN1}). We have to show that
$\sigma=C_*^{-1}$ where $C_*$ is defined by (\ref{constGN2}). Let $(f_n)_{n\in\N}\subset H^1(\R^N)$ be a minimizing sequence. Let
$$
\mu_n=\frac{\|f_n\|_{L^2}^\frac{N-2}{2}}{\|\nabla f_n\|_{L^2}^\frac{N}{2}},\; \lambda_n=\frac{\|f_n\|_{L^2}}{\|\nabla
f_n\|_{L^2}} \;\mbox{ and }\; \forall x\in\R^N,\; v_n(x)=\mu_nf_n(\lambda_nx).
$$
Then $\|v_n\|_{L^2}=\|\nabla v_n\|_{L^2}=1$ and
$J(f_n)=J(v_n)=\|v_n\|_{L^{\alpha+2}}^{-(\alpha+2)}\xrightarrow{n\to\infty}\sigma.$ Let $v_n^*$ be the symmetrization of Schwarz of $|v_n|$ (see Bandle \cite{MR81e:35095}; Berestycki and Lions \cite{MR84h:35054a}, Appendix A.III). Then
$J(v_n^*)\xrightarrow{n\to\infty}\sigma$ and by compactness, $v_{n_\ell}^*\weak v$ as $\ell\tends\infty$ in $H^1_\w(\R^N)$ (and in particular, in $L^{\alpha+2}_\w(\Omega)$ for every subset $\Omega\subset\R^N)$ and $v_{n_\ell}^*\xrightarrow[\ell\to\infty]{L^{\alpha+2}}v,$ for a subsequence $(v_{n_\ell}^*)_\ell\subset(v_n^*)_n$ and for some $v\in H^1_\rad(\R^N).$ Indeed, since $(v_n^*)_{n\in\N}$ is bounded in $H^1_\rad(\R^N)$ and nonincreasing with respect to $|x|,$ then $\forall\ell\in\N$ and $\forall x\in\R^N,$ $|v_{n_\ell}^*(x)|\le C|x|^{-\frac{N}{2}},$ where $C>0$ does not depend on $\ell$ and $x$ (Berestycki and Lions \cite{MR84h:35054a}, Appendix A.II, Radial Lemma A.IV). From this and H{\"o}lder's inequality, we deduce that $\forall\ell\in\N$ and $\forall R>0,$ $\|v_{n_\ell}^*\|_{L^{\alpha+2}(\R^N\setminus B(0,R))}\le CR^{-\frac{N\alpha}{2(\alpha+2)}},$ for a constant $C>0$ which does not depend on $\ell.$ Then $\forall R>0,$ $\|v\|_{L^{\alpha+2}(\R^N\setminus B(0,R))}\le \vliminf_{\ell\to\infty}\|v_{n_\ell}^*\|_{L^{\alpha+2}(\R^N\setminus B(0,R))}\le CR^{-\frac{N\alpha}{2(\alpha+2)}}.$ The strong convergence in $L^{\alpha+2}(\R^N)$ follows easily from the two above estimates and from the compact embedding $H^1(B(0,R))\inj L^{\alpha+2}(B(0,R)),$ which holds for every $R>0.$ Since $\|v_n\|_{L^{\alpha+2}}=\|v_n^*\|_{L^{\alpha+2}},$ it follows that $\|v\|_{L^{\alpha+2}}^{\alpha+2}=\sigma^{-1}$ and then $v\not\equiv0.$ Thus, $J(v)=\sigma$ and $\|v\|_{L^2}=\|\nabla v\|_{L^2}=1.$ It follows that $\forall w\in H^1(\R^N),$ $\frac{d}{dt}J(v+tw)_{|t=0}=0.$ So $v$ satisfies
$-\Delta v+\frac{4-\alpha(N-2)}{N\alpha}v=\sigma\frac{2(\alpha+2)}{N\alpha}|v|^\alpha v,$ in $\R^N.$ Set
$a=\left(\frac{N\alpha}{4-\alpha(N-2)}\right)^\frac{1}{2},$
$b=\left(\frac{2\sigma(\alpha+2)}{4-\alpha(N-2)}\right)^\frac{1}{\alpha}$ and $\forall x\in\R^N,$ $u(x)=bv(ax).$ Then $u\in
H^1_\rad(\R^N)$ is a solution of (\ref{groundstate}) and $J(u)=\sigma.$ By (\ref{groundstate1})--(\ref{groundstate2}), we obtain $J(u)=C_*^{-1}\frac{\|u\|_{L^2}^\alpha}{\|R\|_{L^2}^\alpha}=\sigma$ and $J(R)=C_*^{-1}\ge\sigma$ (since $R$ also satisfies (\ref{groundstate})). Then $\|u\|_{L^2}\le\|R\|_{L^2}$ and so with (\ref{groundstate4}), $\|u\|_{L^2}=\|R\|_{L^2}.$ Hence the result.
\end{vproof}

\section{Necessary condition and sufficient condition for global existence}
\label{cns}

\begin{thm}
\label{thmexplosion}
Let $\lambda >0,$ $\dfrac{4}{N}<\alpha<\dfrac{4}{N-2}$ $(4<\alpha<\infty$ if $N=1)$ and $R$ be the ground state solution of
$(\ref{groundstate}).$ We define for every $a>0,$
\begin{align}
\label{thmexplosion0}
     r_*(a)=\: & \left(\frac{N\alpha}{4-\alpha(N-2)}\right)^\frac{N\alpha-4}{2(4-\alpha(N-2))}
\left(\lambda^{-\frac{1}{\alpha}}\|R\|_{L^2}\right)^\frac{2\alpha}{4-\alpha(N-2)}a^{-\frac{N\alpha-4}{4-\alpha(N-2)}},
\medskip \\
\label{thmexplosion0*}
\gamma_*(a)=\: & \left(\frac{N\alpha-4}{N\alpha}\right)^\frac{N\alpha-4}{2(4-\alpha(N-2))}r_*(a).
\end{align}
\begin{enumerate}
\item
\label{thmexplosion1}
If $\vphi\in H^1(\R^N)$ satisfies
\begin{gather}
\label{thmexplosion11}
\|\vphi\|_{L^2}\le\gamma_*(\|\nabla \vphi\|_{L^2}),
\end{gather}
then the corresponding solution $u\in C((-T_*,T^*);H^1(\R^N))$ of $(\ref{nls})$ is global in time, that is $T^*=T_*=\infty,$
and the following estimates hold.
$$
\begin{array}{rl}
 \forall t\in\R, & \left\{
  \begin{array}{l}
   \|\nabla u(t)\|_{L^2}^2 < \dfrac{2N\alpha}{N\alpha-4}E(\vphi), \medskip \\
   \|\nabla u(t)\|_{L^2} < r_*^{-1}(\|\vphi\|_{L^2}),
  \end{array}
 \right.
\end{array}
$$
where $r_*^{-1}$ is the function defined by $(\ref{rmqthmexplosion11}).$ In particular,
$E(\vphi)>\dfrac{N\alpha-4}{2N\alpha}\|\nabla\vphi\|_{L^2}^2.$
\item
\label{thmexplosion2}
For every $a>0$ and for every $b>0$ satisfying $a>r_*(b),$ there exists $\vphi_{a,b}\in H^1(\R^N)$ with
$\|\vphi_{a,b}\|_{L^2}=a$ and $\|\nabla \vphi_{a,b}\|_{L^2}=b$ such that the associated solution $u_{a,b}\in
C((-T_*,T^*);H^1(\R^N))$ of $(\ref{nls})$ blows up in finite time for both $t>0$ et $t<0.$ In other words, $T^*<\infty$ and
$T_*<\infty.$ Furthermore, $E(\vphi_{a,b})>0 \iff r_*(b)<a<\rho_*(b)$ and $E(\vphi_{a,b})=0 \iff a=\rho_*(b),$ where for every $a>0,$
\begin{gather}
\label{thmexplosion21}
\rho_*(a)=\left(\frac{N\alpha}{4}\right)^\frac{2}{4-\alpha(N-2)}r_*(a).
\end{gather}
Finally, $E(\vphi_{a,b})<\dfrac{N\alpha-4}{2N\alpha}\|\nabla\vphi_{a,b}\|_{L^2}^2.$
\end{enumerate}
\end{thm}

\begin{rmk}
\label{rmqthmexplosion0.7}
Let $\gamma_*$ be the function defined by (\ref{thmexplosion0*}). Set
$$
{\cal A}=\{\vphi\in H^1(\R^N);\; \|\vphi\|_{L^2}\le\gamma_*(\|\nabla \vphi\|_{L^2})\}.
$$
By Theorem \ref{thmexplosion}, for every $\vphi\in{\cal A},$ the corresponding solution of
(\ref{nls}) is global in time and uniformly bounded in $H^1(\R^N).$ It is interesting to note that ${\cal A}$ is an unbounded
subset of $H^1(\R^N).$ So Theorem \ref{thmexplosion} gives a general result for global existence for which we can take initial
values with the $H^1(\R^N)$ norm large as we want.
\end{rmk}

\begin{rmk}
\label{rmqthmexplosion1}
Let $\gamma_*,$ $r_*,$ and $\rho_*$ be the functions defined respectively by (\ref{thmexplosion0*}), (\ref{thmexplosion0}) and
(\ref{thmexplosion21}). It is clear that since $\alpha>\frac{4}{N},$ $\gamma_*,$ $\gamma_*^{-1},$ $r_*,$ $r_*^{-1},$ $\rho_*$
and
$\rho_*^{-1}$ are decreasing and bijective functions from $(0,\infty)$ to $(0,\infty)$ and for every $a>0,$
\begin{align}
\nonumber
\gamma_*^{-1}(a)=\: & \left(\frac{N\alpha-4}{N\alpha}\right)^\frac{1}{2}r_*^{-1}(a), \medskip \\
\label{rmqthmexplosion11}
     r_*^{-1}(a)=\: & \left(\frac{N\alpha}{4-\alpha(N-2)}\right)^\frac{1}{2} 
\left(\lambda^{-\frac{1}{\alpha}}\|R\|_{L^2}\right)^\frac{2\alpha}{N\alpha-4}a^{-\frac{4-\alpha(N-2)}{N\alpha-4}}, \medskip \\
\nonumber
  \rho_*^{-1}(a)=\: & \left(\frac{N\alpha}{4}\right)^\frac{2}{N\alpha-4}r_*^{-1}(a).
\end{align}
So the condition condition (\ref{thmexplosion11}) is equivalent to the condition
$\|\nabla\vphi\|_{L^2}\le\gamma_*^{-1}(\|\vphi\|_{L^2}).$ Furthermore, $\gamma_*<r_*<\rho_*$ and
$\gamma_*^{-1}<r_*^{-1}<\rho_*^{-1}.$
\end{rmk}

\begin{rmk}
\label{rmqthmexplosion2}
Let $\gamma_*,$ $r_*,$ and $\rho_*$ be the functions defined respectively by (\ref{thmexplosion0*}), (\ref{thmexplosion0}) and
(\ref{thmexplosion21}). Then $\gamma_*\xrightarrow{\alpha\searrow\frac{4}{N}}\lambda^{-\frac{1}{\alpha}}\|R\|_{L^2}$ and
$r_*\xrightarrow{\alpha\searrow\frac{4}{N}}\lambda^{-\frac{1}{\alpha}}\|R\|_{L^2}$ (and even,
$\rho_*\xrightarrow{\alpha\searrow\frac{4}{N}}\lambda^{-\frac{1}{\alpha}}\|R\|_{L^2}).$ So we obtain the sharp condition for
global existence, $\|\vphi\|_{L^2}<\lambda^{-\frac{1}{\alpha}}\|R\|_{L^2}$ which coincide with the results obtained by Weinstein
\cite{MR84d:35140,MR87i:35026}. However, we do not know if $\gamma_*$ or $r_*$ are optimum.
\end{rmk}

\section{Proof of Theorem \ref{thmexplosion}}
\label{proofthmexplosion}

In order to prove the blowing up result (\ref{thmexplosion2} of Theorem \ref{thmexplosion}), we need of several lemmas. We
follow the method of Berestycki and Cazenave \cite{MR84f:35120} (see also Cazenave \cite{MR85j:35165} and Cazenave \cite{caz1},
Section 8.2). {\it A priori}, we would expect to use Theorem \ref{explosionX}, that is to construct initial values in $X$ with
nonpositive energy, which is the case for $\alpha=\dfrac{4}{N}.$ But it will not be enough because we have to make blow-up
some solutions whose the initial values have a positive energy.

We define the following functionals and sets. Let $\lambda>0,$ $\omega>0,$ $\beta>0,$ $0<\alpha<\dfrac{4}{N-2}$
$(0<\alpha<\infty$ if $N=1)$ and $\psi\in H^1(\R^N).$
$$
\begin{array}{rl}
 & \left\{
    \begin{array}{l}
      \beta^*(\psi)^\frac{N\alpha-4}{2}=\dfrac{2(\alpha+2)}{\lambda N\alpha}\dfrac{\|\nabla\psi\|_{L^2}^2}
                                          {\|\psi\|_{L^{\alpha+2}}^{\alpha+2}}, \mbox{ if } \psi\not\equiv0, \medskip \\
      Q(\psi)=\|\nabla\psi\|_{L^2}^2-\dfrac{\lambda N\alpha}{2(\alpha+2)}\|\psi\|_{L^{\alpha+2}}^{\alpha+2}, \medskip \\
      S(\psi)=\dfrac{1}{2}\|\nabla\psi\|_{L^2}^2-\dfrac{\lambda}{\alpha+2}\|\psi\|_{L^{\alpha+2}}^{\alpha+2}
                                                                         +\dfrac{\omega}{2}\|\psi\|_{L^2}^2, \medskip \\
      {\cal P}(\beta,\psi)(x)=\beta^\frac{N}{2}\psi(\beta x), \mbox{ for almost every } x\in\R^N,            \medskip \\
      M=\left\{\psi\in H^1(\R^N);\; \psi\not\equiv0\; \mbox{ and }\: Q(\psi)=0\right\},                      \medskip \\
      A=\left\{\psi\in H^1(\R^N);\; \psi\not\equiv0\; \mbox{ and } -\Delta\psi+\omega\psi=\lambda|\psi|^\alpha\psi,
                                                                                   \mbox{ in } \R^N\right\}, \medskip \\
      G=\left\{\psi\in A;\; \forall\phi\in A,\; S(\psi)\le S(\phi)\right\}.
 \end{array}
 \right.
\end{array}
\medskip
$$
Note that by the discussion at the beginning of Section \ref {sharpestimate} and (\ref{groundstate1})--(\ref{groundstate4}),
$M\not=\emptyset,$ $A\not=\emptyset$ and $G\not=\emptyset.$

\begin{lem}
\label{lemthmexplosion0}
We have the following results.
\begin{enumerate}
\item
\label{lemthmexplosion01}
$
\forall\beta>0,\; \beta\not=\beta^*(\psi),\; S({\cal P}(\beta,\psi))<S({\cal P}(\beta^*(\psi),\psi)).
$
\item
\label{lemthmexplosion02}
The following equivalence holds.
$$
\begin{array}{rcl}
 \psi\in G & \iff &
  \left\{
   \begin{array}{l}
    \psi\in M, \\
    S(\psi)=\min\limits_{\phi\in M}S(\phi),
   \end{array}
  \right.
\end{array}
$$
\item
\label{lemthmexplosion03}
Let $m\stackrel{def}{=}\min\limits_{\phi\in M}S(\phi).$ Then
$
\forall\phi\in H^1(\R^N)\;\; with\;\; Q(\phi)<0,\;  Q(\phi)\le S(\phi)-m.
$
\end{enumerate}
\end{lem}

See Cazenave \cite{caz1}, Lemma 8.2.5 for the proof of \ref{lemthmexplosion01}; Proposition 8.2.4 for the proof of
\ref{lemthmexplosion02}; Corollary 8.2.6 for the proof of \ref{lemthmexplosion03}. There is a mistake in the formula (8.2.4) of
this reference. Replace the expression $\lambda^*(u)^\frac{n\alpha-4}{2}=\frac{\alpha+2}{2}\left(\int_{\R^n}|\nabla
u|^2\right)\left(\int_{\R^n}|u|^{\alpha+2}\right)^{-1}$ with $\lambda^*(u)^\frac{n\alpha-4}{2}=\frac{2(\alpha+2)}{n\alpha}
\left(\int_{\R^n}|\nabla u|^2\right)\left(\int_{\R^n}|u|^{\alpha+2}\right)^{-1}.$ \medskip \\
\indent
The proof of \ref{thmexplosion1} of Theorem \ref{thmexplosion} relies on the following lemma.

\begin{lem}
\label{lemthmexplosion3}
Let $I\subseteq\R,$ be an open interval, $t_0\in I,$ $p>1,$ $a>0,$ $b>0$ and $\Phi\in C(I;\R^+).$ We set, $\forall
x\ge0,$ $f(x)=a-x+bx^p,$ $\udl x=(bp)^{-\frac{1}{p-1}}$ and $b_*=\dfrac{p-1}{p}\udl x.$ Assume that $\Phi(t_0)<\udl x,$ $a\le
b_*$ and that $f\circ\Phi>0.$ Then, $\forall t\in I,$ $\Phi(t)<\udl x.$
\end{lem}

\begin{proof*}
Since $\Phi(t_0)<\udl x$ and $\Phi$ is a continuous function, there exists $\eta>0$ with $(t_0-\eta,t_0+\eta)\subseteq I$
such that, $\forall t\in(t_0-\eta,t_0+\eta),$ $\Phi(t)<\udl x.$ If $\Phi(t_*)=\udl x$ for some $t_*\in I,$ then
$f\circ\Phi(t_*)=f(\udl x)=a-b_*\le0.$ But $f\circ\Phi>0.$ Then, $\forall t\in I,$
$\Phi(t)<\udl x.$
\medskip
\end{proof*}

The proof of \ref{thmexplosion2} of Theorem \ref{thmexplosion} makes use the following lemma.

\begin{lem}
\label{lemthmexplosion4}
Let $\lambda>0,$ $\omega>0$ and $\dfrac{4}{N}<\alpha<\dfrac{4}{N-2}$ $(4<\alpha<\infty$ if $N=1).$ We set for every $\beta>0$
and for every $\psi\in H^1(\R^N),$ $\vphi_\beta={\cal P}(\beta,\psi).$ Let $u_\beta\in C((-T_*,T^*);H^1(\R^N))$ be the solution
of $(\ref{nls})$ with initial value $\vphi_\beta.$ Then we have, $\forall\psi\in G,$ $\forall\beta>1,$ $T_*<\infty$ and
$T^*<\infty.$
\end{lem}

\begin{proof*}
Let $\psi\in G.$ By (\ref{groundstate3}), we have
\begin{gather*}
\beta^*(\vphi_\beta)^\frac{N\alpha-4}{2}=\dfrac{2(\alpha+2)}{\lambda
N\alpha}\frac{\beta^2\|\nabla\psi\|_{L^2}^2}{\beta^\frac{N\alpha}{2}\|\psi\|_{L^{\alpha+2}}^{\alpha+2}}, \\
Q(\vphi_\beta)=\beta^2\left(\|\nabla\psi\|_{L^2}^2-\frac{\lambda
N\alpha}{2(\alpha+2)}\beta^\frac{N\alpha-4}{2}\|\psi\|_{L^{\alpha+2}}^{\alpha+2}\right), \\
\frac{\lambda N\alpha}{2(\alpha+2)}\|\psi\|_{L^{\alpha+2}}^{\alpha+2}=\|\nabla\psi\|_{L^2}^2.
\end{gather*}
So, $\beta^*(\vphi_\beta)^\frac{N\alpha-4}{2}=\beta^{-\frac{N\alpha-4}{2}},$
$Q(\vphi_\beta)=-\beta^2\|\nabla\psi\|_{L^2}^2\left(\beta^\frac{N\alpha-4}{2}-1\right)$ and $\beta^*(\psi)=1.$ From these three
last equalities, from \ref{lemthmexplosion01} and \ref{lemthmexplosion02} of Lemmas \ref{lemthmexplosion0} and by conservation of
charge and energy, we have
\begin{gather}
\label{demothmexplosion22}
\forall\beta>1,\; Q(\vphi_\beta)<0, \\
\label{demothmexplosion23}
\forall\beta\not=1,\; S(\vphi_\beta)<S(\psi)\equiv m, \\
\label{demothmexplosion24}
\forall\beta>0,\; \forall t\in(-T_*,T^*),\; S(u_\beta(t))=S(\vphi_\beta).
\end{gather}
By continuity of $u_\beta,$ by (\ref{demothmexplosion22})--(\ref{demothmexplosion24}) and from \ref{lemthmexplosion03} of Lemma
\ref{lemthmexplosion0}, we have for every $\beta>1,$
\begin{gather}
\label{demothmexplosion25}
\forall t\in(-T_*,T^*),\; Q(u_\beta(t))\le S(\vphi_\beta)-m<0.
\end{gather}
Set $\forall t\in(-T_*,T^*),$ $h(t)=\|xu_\beta(t)\|_{L^2}^2.$ Then we have by Glassey \cite{MR57:842}, $h\in C^2((-T_*,T^*);\R)$
and $\forall t\in(-T_*,T^*),$ $h''(t)=8\|\nabla u_\beta(t)\|_{L^2}^2-\frac{4\lambda
N\alpha}{\alpha+2}\|u_\beta(t)\|_{L^{\alpha+2}}^{\alpha+2}\equiv8Q(u_\beta(t)).$ So with (\ref{demothmexplosion25}),
$$
\forall t\in(-T_*,T^*),\; h''(t)\le8(S(\vphi_\beta)-m)<0,
$$
for every $\beta>1.$ It follows that $T_*<\infty$ and $T^*<\infty.$ Hence the result.
\medskip
\end{proof*}

\begin{vproof}{of Theorem \ref{thmexplosion}.}
We proceed in two steps. \\
{\bf Step 1.} We have \ref{thmexplosion1}. \\
Let $C_*$ be the constant defined by (\ref{constGN2}). We set $I=(-T_*,T^*),$ $t_0=0,$ $p=\frac{N\alpha}{4},$
$a=\|\nabla\vphi\|_{L^2}^2,$ $b=\frac{2\lambda}{\alpha+2}C_*\|\vphi\|_{L^2}^\frac{4-\alpha(N-2)}{2},$
$\udl x=(bp)^{-\frac{1}{p-1}},$ $b_*=\frac{p-1}{p}\udl x,$ $\forall t\in I,$ $\Phi(t)=\|\nabla u(t)\|_{L^2}^2$ and for any $x\ge0,$ $f(x)=a-x+bx^p.$ Then by conservation of energy, by Proposition \ref{constGN} and by conservation of charge, we have
$$
\begin{array}{rcl}
\forall t\in I,\; \|\nabla u(t)\|_{L^2}^2 & = & 2E(\vphi)+\dfrac{2\lambda}{\alpha+2}\|u(t)\|_{L^{\alpha+2}}^{\alpha+2}\medskip\\
                                          & < & \|\nabla\vphi\|_{L^2}^2+\dfrac{2\lambda}{\alpha+2}C_*\|\vphi\|_{L^2}
                                                ^\frac{4-\alpha(N-2)}{2} (\|\nabla u(t)\|_{L^2}^2)^\frac{N\alpha}{4}.
\end{array}
$$
And so, $\forall t\in I,$ $a-\|\nabla u(t)\|_{L^2}^2+b(\|\nabla u(t)\|_{L^2}^2)^p>0,$
that is $f\circ\Phi>0.$ Furthermore, $\Phi(t_0)\equiv a\le b_*<\udl x.$ Indeed, by Remark \ref{rmqthmexplosion1}, we have
$$
\Phi(t_0)\le b_*\iff\|\nabla\vphi\|_{L^2}\le\gamma_*^{-1}(\|\vphi\|_{L^2})\iff\|\vphi\|_{L^2}\le\gamma_*(\|\nabla\vphi\|_{L^2}).
$$
So by Lemma \ref{lemthmexplosion3}, $\Phi(t)<\udl x\equiv[r_*^{-1}(\|\vphi\|_{L^2})]^2,$ $\forall t\in I.$ Thus, $I=\R$ and for
every $t\in\R,$
$$
\|\nabla u(t)\|_{L^2}<r_*^{-1}(\|\vphi\|_{L^2}).
$$
It follows from conservation of charge and energy, (\ref{constGN1}), (\ref{constGN2}), and the above inequality, that
\begin{align*}
  \forall t\in\R,\; E(\vphi)
& \ge \frac{1}{2}\left(\|\nabla u(t)\|_{L^2}^2-\frac{2\lambda}{\alpha+2}C_*\|\vphi\|_{L^2}^\frac{4-\alpha(N-2)}{2}
                                                                        \|\nabla u(t)\|_{L^2}^\frac{N\alpha}{2}\right) \\
&  =  \frac{1}{2}\|\nabla u(t)\|_{L^2}^2\left(1-\frac{4}{N\alpha}\left[r_*^{-1}(\|\vphi\|_{L^2})\|\nabla
                                                       u(t)\|_{L^2}^{-1}\right]^{-\frac{N\alpha-4}{2}}\right) \medskip \\
&  >  \frac{1}{2}\|\nabla u(t)\|_{L^2}^2\left(1-\frac{4}{N\alpha}\right) \medskip \\
&  =  \frac{N\alpha-4}{2N\alpha}\|\nabla u(t)\|_{L^2}^2.
\end{align*}
Hence \ref{thmexplosion1}. \\
{\bf Step 2.} We have \ref{thmexplosion2}. \\
Let $R$ be the ground state solution of (\ref{groundstate}). Let first remark from the assumptions and from Remark
\ref{rmqthmexplosion1}, we have $b>r_*^{-1}(a).$ We set
\begin{center}
$\nu=[r_*^{-1}(a)]^\frac{N}{2}\left(\frac{4-\alpha(N-2)}{N\alpha}\right)^\frac{N}{4}a^{-\frac{N-2}{2}}\|R\|_{L^2}^{-1},$ \quad
$\omega=[r_*^{-1}(a)]^2\frac{4-\alpha(N-2)}{N\alpha}a^{-2}=\left(\lambda^{-\frac{1}{\alpha}}\|R\|_{L^2}a^{-1}\right)
^\frac{4\alpha}{N\alpha-4},$
\end{center}
and for every $x\in\R^N,$ $\psi(x)=\nu R(\sqrt\omega x).$ Then $\psi\in\vS_\rad(\R^N)\cap A.$ Since $R$ satisfies
(\ref{groundstate})--(\ref{groundstate4}), it follows that $\psi\in G.$ Furthermore, $\|\psi\|_{L^2}=a$ and
$\|\nabla\psi\|_{L^2}=r_*^{-1}(a).$ Let $\beta=\frac{b}{r_*^{-1}(a)}>1.$ Set for every $x\in\R^N,$
$\vphi_{a,b}(x)=\vphi_\beta(x)={\cal P}(\beta,\psi)(x).$ In particular, $\vphi_{a,b}\in\vS_\rad(\R^N)$ and $\vphi_{a,b}$
satisfies
$$
-\Delta\vphi_{a,b}+\omega\beta^2\vphi_{a,b}=\lambda\beta^{-\frac{N\alpha-4}{2}}|\vphi_{a,b}|^\alpha\vphi_{a,b}, \mbox{ in } \R^N.
$$
Denote $u_{a,b}\in C((-T_*,T^*);H^2(\R^N)\cap X_\rad)$ the solution of $(\ref{nls})$ with initial value $\vphi_{a,b}.$ Then by
Lemma \ref{lemthmexplosion4}, $T_*<\infty$ and $T^*<\infty.$ Moreover, $\|\vphi_{a,b}\|_{L^2}=a,$
$\|\nabla\vphi_{a,b}\|_{L^2}=b$ and by (\ref{groundstate3}),
$$
\begin{array}{rcl}
E(\vphi_{a,b}) & = & \dfrac{1}{2}\|\nabla\vphi_{a,b}\|_{L^2}^2-\dfrac{\lambda}{\alpha+2}\|\vphi_{a,b}\|_{L^{\alpha+2}}^{\alpha+2}
                                                                                                                    \medskip \\
               & = & \dfrac{1}{2}\|\nabla\vphi_{a,b}\|_{L^2}^2-\dfrac{\lambda}{\alpha+2}\beta^\frac{N\alpha}{2}
                                                                                 \|\psi\|_{L^{\alpha+2}}^{\alpha+2} \medskip \\
               & = & \dfrac{\|\nabla\vphi_{a,b}\|_{L^2}^2}{2N\alpha}\left(N\alpha-4\beta^\frac{N\alpha-4}{2}\right) \medskip \\
               & = & \dfrac{\|\nabla\vphi_{a,b}\|_{L^2}^2}{2N\alpha}\left(N\alpha-4\beta^\frac{N\alpha-4}{2}\right).
\end{array}
$$
By Remark \ref{rmqthmexplosion1}, it follows that
$$
E(\vphi_{a,b})\le0\iff\beta\ge\left(\frac{N\alpha}{4}\right)^\frac{2}{N\alpha-4}\iff
b\ge\left(\frac{N\alpha}{4}\right)^\frac{2}{N\alpha-4}r_*^{-1}(a)\equiv\rho_*^{-1}(a)\iff a\ge\rho_*(b).
$$
Hence the result.
\bigskip
\end{vproof}

\baselineskip .5cm

\noindent
{\large\bf Acknowledgments} \\
The author would like to thank his thesis adviser, Professor Thierry Cazenave, for his suggestions and encouragement.

\baselineskip .0cm

\bibliographystyle{abbrv}
\bibliography{BiblioPaper2}

\end{document}